\input Tex-document.sty

\def\firstpage{373}

\pageno=373

\def\shorttitle{Dynamics in Two Complex Dimensions}
\def\shortauthor{J. Smillie}

\title{\centerline{Dynamics in Two}
\centerline{Complex Dimensions}}

\author{J. Smillie\footnote{\eightrm $^\ast$}{\eightrm Department of Mathematics, Malott Hall, Cornell
University, Ithaca, NY, USA. E-mail: smillie@math.cornell.edu}}

\vskip 7mm

\centerline{\boldnormal Abstract}

\vskip 4.5mm

{\narrower \ninepoint \smallskip We describe results on  the
dynamics of polynomial diffeomorphisms of
  ${\bf C^2}$ and draw connections with the dynamics of polynomial maps of
${\bf C}$ and the dynamics of polynomial diffeomorphisms of ${\bf
R^2}$ such as the H\'enon family.

\vskip 4.5mm



}

\vskip 10mm

\head{1. Introduction}

     The  subject of this article is part of the larger subject
area of higher dimensional complex dynamics. This larger area
includes the dynamical study of holomorphic maps of complex
projective space, automorphisms of $K3$ surfaces, birational maps,
automorphisms of ${\bf C^n}$ and higher dimensional Newton's
method. Our particular topic of research is polynomial
automorphisms of ${\bf C^2}$. This area is particularly
interesting because of its connections to some fundamental
questions of dynamical systems via two real dimensional dynamics
and because of its connection to some powerful techniques via one
dimensional complex dynamics. I will begin by describing some of
these connections. The reader is encouraged to consult [17] for a
more thorough discussion of the historical background summarized
here.

   Over one hundred years after Poincar\'e observed chaotic behavior in the
dynamics of surface diffeomorphisms the problem of creating a
comprehensive theory of the dynamics of diffeomorphisms remains
unsolved. Though the objective is to create a theory that would
apply to diffeomorphisms in any dimension the focus remains on the
two dimensional case. On the one hand the chaotic behavior which
makes these problems challenging first appears for diffeomorphisms
in dimension two, on the other hand there is a sense that if the
tools can be developed to solve the problem in dimension two then
the higher dimensional problem will be approachable. There are
reasons to believe that if the tools can be developed to
thoroughly analyze one specific interesting family of
diffeomorphisms then one would be in a good position to attack the
general problem. If we were to suggest a family to play the role
of a ``test case'' there is one particular family which stands
out.  This is the family of diffeomorphisms of ${\bf R^2}$
introduced by the French astronomer H\'enon:
$$f_{a,b}(x,y)=(a-by-x^2,x).$$

The parameter $b$ is the Jacobian determinant of $f_{a,b}$. When
$b\ne 0$ these maps are diffeomorphisms. When $b=0$ then $f_{a,0}$
is a map with a one dimension range and the behavior of $f_{a,0}$
is essentially that of the quadratic unimodal map $f_a(x)=a-x^2$.

In singling out the H\'enon family we are following a well
established tradition. This family has appeared often both in the
physics and mathematics literature. It has been studied
theoretically and numerically.

Virtually all interesting dynamical behavior which is known to
occur for two dimensional diffeomorphisms is known to occur in
this family. H\'enon's original question involved an apparent
strange attractor, and this is the first family in which the
existence of strange attractors was proved ([11]). For certain
parameter values this family exhibits hyperbolic behavior such as
the Smale horseshoe ([14]). For other parameters it exhibits
persistently nonhyperbolic behavior ([19]).

There is also a great deal that is not understood about the
H\'enon family. Despite the fact that many different types of
dynamic behavior occur it is not known whether the union of these
behaviors accounts for a large set of parameter values. There are
also open questions about how the complexity of behavior varies
with the parameters. When $a\ll0$ the behavior is non-chaotic.
When $a\gg0$, $f_{a,b}$ exhibits a horseshoe, a model for chaotic
behavior. What happens for intermediate values? How is chaos
created? (cf [13])

Another reason for looking at the H\'enon family is its connection
with the one dimensional family of unimodal maps $f_a$. One
dimensional diffeomorphisms exhibit only regular behavior but one
dimensional maps exhibit a wealth of chaotic behavior. In contrast
to the situation for the H\'enon family, the most fundamental
questions for the unimodal family have been answered. In the
language of [17] the family $f_a$ provides a ``qualitatively
solvable model of chaos" which is to say that there is a good
understanding of attractors, strange and otherwise, for large sets
of parameters and there is a good understanding of the transition
to chaos.

The quadratic family is distinguished in the family of unimodal
maps because it has a natural extension to the complex numbers. In
the family $f_a(x)=a-x^2$ both $x$ and $a$ can be taken to be
complex. The use of complex methods stands out as a reason for the
success of the analysis of the quadratic family and unimodal maps
more generally. While there are important results about unimodal
maps that do not use complex techniques, these techniques do play
a central role in the monotonicity results and in the analysis of
attractors for the quadratic family.

Because the H\'enon family is also given by polynomial equations
it also has a natural complex extension.  My first introduction to
the importance of the complex H\'enon family was through lectures
of J. H. Hubbard in the mid 1980's. Another contributor who
brought  new ideas to the subject was N. Sibony. Hubbard and his
co-authors as well as Fornaess and Sibony and many others have
continued to make fundamental contributions to this area and it is
not possible to do justice to all of this work in the space
provided.  I will focus here on work that was carried out jointly
with E. Bedford and, in some cases, M. Lyubich over the past 15
years.

\head{2. Basic definitions in one and two variables}

The fundamental paper of Friedland and Milnor [15] shows that a
natural class of holomorphic diffeomorphisms to consider is the
family of polynomial diffeomorphisms of ${\bf C^2}$. This class
contains the H\'enon family and the tools that we use to analyze
the H\'enon family work equally well for all diffeomorphisms in
this class. In studying polynomial maps of ${\bf C}$ one focuses
on those of degree greater than one because these exhibit chaotic
behavior.  One way of quantifying chaotic behavior is through the
topological entropy, $h_{top}(f)$. In one complex dimension the
entropy is the logarithm of the degree so the distinction between
degree one and higher degree is the distinction between entropy
zero and positive entropy.

For polynomial diffeomorphisms in dimension two the algebraic
degree is not a conjugacy invariant and hence not a dynamical
invariant.  One way to create a conjugacy invariant is to define
the following ``dynamical degree":
$$d=\lim_{n\to\infty} ({\rm algebraic\ degree\ } f^n)^{1\over n}.$$

It is again true that the topological entropy of a complex
diffeomorphism is the logarithm of its dynamical degree, so
dynamical degree seems to be the appropriate two dimensional
analog of degree. The H\'enon diffeomorphisms have the property
that the algebraic degree of $f^n$ is $2^n$ so the dynamical
degree is two. Friedland and Milnor show that any diffeomorphism
with dynamical degree one is conjugate to an affine or elementary
diffeomorphism. They also show that a diffeomorphism with
dynamical degree greater than one is, like the H\'enon
diffeomorphism, conjugate to an explicit diffeomorphism whose
actual degree is equal to its dynamical degree.   When we refer to
the degree of a diffeomorphism we will mean the dynamical degree.
We make the standing assumption that all of our polynomial
diffeomorphisms have degree greater than one.

Let us review some standard definitions for polynomial maps.  Let
$f:{\bf C}\to {\bf C}$ be a polynomial map with degree $d>1$. The
set $K$ is the set of points with bounded orbits. The Julia set,
$J$ is the boundary of $K$. In dimension one all recurrent
behavior is contained in $K$. All chaotic recurrent behavior is
contained in $J$. The ease with which this set can be defined
leaves one unprepared for the range of intricate behavior that it
exhibits.

Let $f:{\bf C^2}\to{\bf C^2}$ be a polynomial diffeomorphism with
dynamical degree $d>1$. The set $K^+$ is the set of points with
bounded forward orbits. Following Hubbard we take the set $K^-$ to
be the set of points with bounded backward orbits. The sets
$J^\pm$ are defined as the boundaries of $K^\pm$. The set $J$ is
defined to be $J^+\cap J^-$. In dimension two all chaotic
recurrent behavior is contained in $J$. Thus $J$ seems to be a
good analog of the one dimensional Julia set. (In fact there is an
alternative analog of $J$ but we will not deal with that here.)

Let $p$ be a periodic saddle point of period $n$ in ${\bf C^2}$.
Let $W^u_p$ denote the unstable manifold of $p$. This is the set
of points that converge to $p$ under iteration of $f^{-1}$. Since
this definition involves $f^{-1}$ it is less clear what the one
variable analog should be. Let us examine the situation more
carefully. The set $W^u_p$ is holomorphically equivalent to ${\bf
C}$. We can find a parameterization $\phi_p:{\bf C}\to W^u_p$
which satisfies the functional equation
$f^n(\phi_p(z))=\lambda\cdot z$ where $\lambda$ is the expanding
eigenvalue of $Df^n_p$. Now if $p$ is a periodic point in ${\bf
C}$ then the functional equation still makes sense. A function
$\phi_p$ which satisfies this equation is called a linearizing
coordinate and this is a good analog of the parameterized unstable
manifold in two dimensions.

Hubbard made the key observation that this construction gives a
natural way to draw pictures of the sets $W^u_p\cap K^+$ in two
variables and a natural way to compare them to the pictures of $K$
in one variable. In both cases we identify a region in ${\bf C}$
with the computer screen and choose a color scheme where the color
for a pixel corresponding to $z$ is related to the rate of escape
of $\phi_p(z)$. The general convention is that points that do not
escape (those points in $\phi_p^{-1}(K)$) are colored black. (See
[http://www.math.cornell.edu/\~{ }dynamics/].)

There is an abstract construction which makes it easier to compare
invertible systems such as diffeomorphisms with non-invertible
systems. Given a non-invertible system such as $f: {\bf C}\to {\bf
C}$ there is a closely related invertible system called the
natural extension.  Let us denote this by ${\hat f}:{\bf \hat
C}\to {\bf\hat C}$. The points in ${\bf\hat C}$ consist of
sequences $(\ldots z_{-1}, z_0, z_1\ldots)$ such that
$f(z_j)=z_{j+1}$. The map ${\hat f}$ acts by shifting such a
sequence to the left.

The natural extension gives us a way of justifying the analogy
between linearizing coordinates and unstable manifolds.
Corresponding to a periodic saddle point $p$ in ${\bf C}$ there is
a unique periodic point $\hat p$ in ${\bf\hat C}$. Since $\hat f$
is invertible we can make sense of the unstable manifold
$W^u_{\hat p}$ and the linearizing coordinate can be used to
parameterize this unstable manifold.

Though ${\bf\hat C}$ contains ``leaves" such as $W^u_{\hat p}$ it
is a mistake to think of ${\bf\hat C}$ as a lamination. When $f$
is expanding ${\bf\hat C}$ is a lamination near the Julia set but
the more complicated the dynamics of $f$, the more degenerate this
structure becomes. This complexity arises from recurrent behavior
of the critical point for $f$. This suggests a certain connection
between regularity of unstable manifolds in two dimensions and
recurrence of critical points in one dimension that we will return
to later.

Since  points in ${\bf \hat C}$ have bounded backward orbits, we
should think of ${\bf\hat C}$ as an analog of $J^-$. Let $f_a$ be
an expanding one dimensional map and consider a diffeomorphism
$f_{a',b}$ with $b$ small and $a'$ close to $a$.  Hubbard and
Oberste-Vorth ([17]) show that $J^-$ is topologically conjugate to
the corresponding ${\bf \hat C}$. When the one dimensional map $f$
is not expanding the relation between ${\bf \hat C}$ and any
particular $J^-$ should be viewed as metaphorical rather than
literal.

\head{3. Potential theory and Pluri-potential theory}

A standard construction in potential theory associates to nice
sets a measure $\mu$ called the harmonic measure or equilibrium
measure. The harmonic measure associated to the Julia set turns
out to be a measure of dynamical interest. The potential theory
construction starts with the Green function. The Green function of
$K$ has a dynamical description:
  $$G(p)=\lim_{n\to\infty}{1\over {d^n}} \log^+|f^n(p)|.$$
The Green function is non-negative and equal to zero precisely on
the set $K$. The harmonic measure $\mu$ is obtained by applying
the Laplacian to $G$. The support of $\mu$ is the boundary of $K$
which is the set $J$.
  The connection between polynomial maps and potential theory first appears in
the work of Brolin ([12]). It reappears in a paper of Manning
([18]) and is nicely summarized in [20].

The harmonic measure has connections to entropy and to the
connectivity of $J$. These connections do not play a major role in
the one dimensional theory because entropy and connectivity can be
approached more directly. In the two dimensional theory these
connections are much more important.

The entropy of the measure $\mu$, $h(\mu)$, happens to be $\log d$
which is equal to the topological entropy of the map. The
topological entropy dominates the measure theoretic entropy of any
invariant measure. A measure for which equality holds is called a
measure of maximal entropy. For polynomial maps of ${\bf C}$ the
measure $\mu$ can be characterized as the unique  measure of
maximal entropy.

The dimension of a measure $\nu$, $\dim_H(\nu)$, is the minimum of
the Hausdorff dimensions of subsets of full $\nu$ measure. The
dimension of the harmonic measure of a planar set is always less
than or equal to one. If the set is connected then the dimension
is one. For Julia sets the converse is true: the dimension of the
harmonic measure is one if and only if $J$ is connected.

The Lyapunov exponent, $\lambda(\mu)$, of $f$ with respect to an
ergodic measure measures the rate of growth of tangent vectors
under iteration (for a set of full $\mu$ measure). The Lyapunov
exponent is related to Hausdorff dimension of the measure by the
formula:
$$\dim_H(\mu)=h(\mu)/\lambda(\mu).$$
Since $h(\mu)$ is $\log d$,  $\lambda(\mu)=\log d$ if and only if
$J$ is connected. We will return to this in the next section.

In dimension two we have two rate of escape functions:
$$G^\pm(p)=\lim_{n\to\infty}{1\over {d^n}} \log^+|f^{\pm n}(p)|.$$
  Potential theory in one
variable centers on the behavior of the Laplacian. The Laplacian
is not holomorphically invariant in two variables but it has a
close relative which is. This is the operator $dd^c$ which takes
real valued functions to real two forms. The $d$ that appears here
is just the exterior derivative and the $d^c$ is a version of the
exterior derivative twisted by using the complex structure. In one
variable we have:
$$dd^cg=(\triangle g )dx\wedge dy.$$
Not only is  $dd^c$ holomorphically natural but it is well defined
on complex manifolds of any dimension. Of course, in the two
variable context as in the one variable, the functions to which
these operators are applied are not smooth and the result has to
be interpreted appropriately. The theory connected with the
operator $dd^c$ is referred to as pluripotential theory. It was an
observation of Sibony that the methods of pluripotential can be
profitably applied to the complex H\'enon diffeomorphisms.

Define $\mu^\pm={1\over{2\pi}}dd^c G^\pm$. These are dynamically
significant currents supported on $J^\pm$.  Define
$\mu=\mu^+\wedge\mu^-$. This measure $\mu$ is the analog of the
harmonic measure in one dimension. The following result suggest
that ``$\mu$" defined above is the dynamical analog as well as the
pluripotential theoretic analog.

{\bf Theorem 3.1 ([4])} \it The measure $\mu$ is the unique
measure of maximal entropy. \rm

\head{4. Connectivity and critical points}

We want to consider the way in which the dynamical behavior of a
polynomial diffeomorphism such as $f_{a,b}$ depends on the
parameter. Looking at pictures of $W^u_p\cap K^+$ shows that there
are indeed many things that do change. If we want to focus on one
fundamental property we might start by looking at connectivity. In
one variable the connectivity of the Julia set of $f_a$ defines
the Mandelbrot set which is the fundamental object of study for
quadratic maps. In two variables there are several notions of
connectivity that we could consider. The following has proved
useful. We say that $f$ is {\it stably/unstably connected} if
$W_p^{s/u}\cap K^{-/+}$ has no compact components. We can ask
about the relation between stable connectivity, unstable
connectivity and the connectivity of $J$. A priori the property of
being stably/unstably connected depends on the saddle point $p$.
In fact we show that these properties are independent of $p$.

Let us look at the situation in one variable. The basic result
about connectivity is the following.

{\bf Theorem 4.1 (Fatou)} \it Let $f$ be a polynomial map of ${\bf
C}$. Then $J$ is connected if and only if every critical point of
$f$ has a bounded orbit. \rm

The following formula makes a connection between the Lyapunov
exponent and critical points ([20]):
$$ \lambda(\mu)=\log d+\sum_{\{c_j:f'(c_j)=0\} } G(c_j).$$
The function $G$ is non-negative and zero precisely on the set
$K$. In light of the theorem above we see that $J$ is connected if
and only if the Lyapunov exponent is $\log d$. This proves an
assertion made in Section 1 about the relation between the
Lyapunov exponents of $\mu$ and the connectivity of $J$.

In two variables there are two Lyapunov exponents,
$\lambda^\pm(\mu)$, of $f$ with respect to harmonic measure. The
following result establishes the connection between
stable/unstable connectivity and these exponents.

{\bf Theorem 4.2 ([7])} \it We have $\lambda^+(\mu)\ge\log d$; and
$\lambda^+(\mu)=\log d$ if and only if $f$ is unstably connected.
Similarly $\lambda^-(\mu)\le-\log d$; and $\lambda^-(\mu)=-\log d$
if and only if $f$ is stably connected. \rm

It is clear from this result that neither exponent is zero. Pesin
theory shows that stable and unstable manifolds exist for $\mu$
almost every point. Let ${\cal C}^u$ be the set of critical points
of the restriction of $G^+$ to these unstable manifolds. We define
${\cal C}^s$ in the corresponding way.

{\bf Theorem 4.3 ([7])} \it The diffeomorphism $f$ is unstably
connected if and only if ${\cal C}^u=\emptyset$. The
diffeomorphism $f$ is stably connected if and only if ${\cal
C}^s=\emptyset$. \rm

In [6] we prove an analog of the critical point formula where the
role of the critical point is played by critical points is played
by ${\cal C}^u$. This formula leads to proofs of the two theorems
above.

The following result makes the connection between stable and
unstable connectivity and the connectivity of $J$. Note that in
this two variable situation the Jacobian of $f$ enters the
picture.

{\bf Theorem 4.4 ([7])} \it If $|\det Df|<1$ then $f$ is never
stably connected. In this case $J$ is connected if and only if $f$
is unstably connected. If $|\det Df|=1$ then $f$ is stably
connected iff $f$ is unstably connected iff $J$ is connected. \rm

(The case $|\det Df|> 1$ is analogous to the case $|\det Df|< 1$.)
The Jacobian enters the proof through the  relation:
$\lambda^+(\mu)+\lambda^-(\mu)=\log |\det Df|$. We see for example
that
  $|\det Df|< 1$ implies that $\lambda^-(\mu)< -\log d$ which, by Theorem 4.2
implies that $f$ is unstably disconnected.

Using this result J. H. Hubbard and K. Papadantonakis have
developed a computer program that uses the set ${\cal C}^u$ to
draw pictures of the connectivity locus in parameter space. (See
[http://www.math.cornell.edu/\~{ }dynamics/].)

\head{5. The boundary of the horseshoe locus}

Hyperbolic behavior, as exhibited by the horseshoe, is
structurally stable. This implies that the set of $(a,b)$ for
which $f_{a,b}$ exhibits a horseshoe is open. Let us call this set
the horseshoe locus. Standard techniques from dynamical systems
can be used to analyze the dynamical behavior inside the horseshoe
locus. These techniques break down on the boundary of the
horseshoe locus however. By contrast complex techniques from [4],
[9] and [10] can be applied on the closure of the horseshoe locus.
Thus the analysis of this boundary provides a setting for
demonstrating that these techniques derived from complex analysis
are not without interest in the real setting.

Let us look at the one dimensional case $f_a$. We say that $f_a$
exhibits a horseshoe if $f_a|J_a$ is expanding and topologically
conjugate to the one sided two shift. The horseshoe locus here is
the set $a>2$. The boundary of the horseshoe locus is $a=2$. The
map $f_2$ is the well known example of Ulam and von Neumann. The
failure of expansion is demonstrated by the fact that the critical
point $0$ is in the Julia set, $[-2,2]$. In fact the critical
point maps to the fixed point $-2$ after two iterates.

The following result describes the failure of hyperbolicity on the
boundary of the horseshoe locus for H\'enon diffeomorphisms. Note
that the property of eventually mapping to the fixed point $p$ in
dimension one corresponds to belonging to $W^s_p$ in dimension
two.

{\bf Theorem 5.1 ([10])}. \it For $f_{a,b}$ on the boundary of the
horseshoe locus there are fixed points $p$ and $q$ so that $W^s_p$
and $W^u_q$ have a quadratic tangency. When $b>0$ we have $p=q$.
When $b<0$, $p\ne q$. \rm

The next result gives additional information about the precise
nature of the dynamics of maps on the boundary of the horseshoe
locus:

{\bf Theorem 5.2 (Bedford-Smillie)} \it For any $(a,b)$ in the
boundary of the horseshoe locus the restriction of $f_{a,b}$ to
its non-wandering set is conjugate to the full two-shift with
precisely two orbits identified. Given $(a,b)$ and $(a',b')$ in
the boundary of the horseshoe locus, the restrictions of $f_{a,b}$
and $f_{a',b'}$ to their non-wandering sets are conjugate if and
only if $b$ and $b'$ have the same sign. \rm

There are many techniques which work only for $b$ small. Note that
that the result above applies for all values of $b$ including the
volume preserving case $b=\pm 1$.

We can ask how the dynamics of $f_{a,b}$ for $(a,b)$ on the
boundary of the horseshoe regions $b>0$ and $b<0$ compares with
the dynamics of $f_2$ which corresponds to the boundary of the
horseshoe region when $b=0$. The sets $J_{a,b}$ for $b\ne0$ are
totally disconnected while the set $J_2$ is connected. In
particular the inverse limit system $\hat J_2$ is not conjugate to
either system with $b\ne 0$. This is an example where the insights
gained from looking at the inverse limit system need to be
interpreted cautiously.

I will touch on the techniques used in the proofs of these
theorems. Our fundamental approach to proving these results was to
exploit the relationship between the real mapping $f_{a,b}:{\bf
R^2}\to{\bf R}^2$ and its complex extension $f_{a,b}:{\bf
C^2}\to{\bf C}^2$. In passing from ${\bf C^2}$ to ${\bf R^2}$
something may be lost. The first question to ask is how much
chaotic behavior do we lose? One way to measure this is through
the topological entropy function. If we denote $f_{a,b}:{\bf
R^2}\to{\bf R}^2$ by $f_{\bf R}$ and $f_{a,b}:{\bf C^2}\to{\bf
C}^2$ by $f_{\bf C}$ then we have
$$h_{top}(f_{\bf R})\le h_{top}(f_{\bf C})=\log 2.$$

If we want to study the real H\'enon diffeomorphisms most closely
connected to their complex extensions we should focus our
attention on those $f$ with $h_{top}(f_{\bf R})=\log 2$. We say
that these examples have {\it maximal entropy}. This is an
interesting set to look at. The horseshoe locus is contained in
the maximal entropy locus but the maximal entropy locus is larger
than the horseshoe locus. The horseshoe locus is open, and the
maximal entropy locus is  closed. In particular the maximal
entropy locus contains the boundary of the horseshoe locus.

For maximal entropy diffeomorphisms the relation between the real
and complex dynamics is as close as one could want:

{\bf Theorem 5.3 ([4])} \it $f_{\bf R}$ has maximal entropy if and
only if $J$ is contained in  ${\bf R^2}$. \rm

This theorem is a consequence of the fact that $\mu$ is the unique
measure of maximal entropy and the fact that the support of $\mu$
is contained in $J$. The fact that the real and complex dynamics
are closely related for this class of maps means that it is a good
starting point for applying complex techniques to the real case.
It also provides us with useful techniques from harmonic analysis.
For example the Green functions of real sets satisfy certain
growth conditions, and these translate into conditions insuring
expansion and regularity of unstable manifolds. This allows us to
show that maximal entropy diffeomorphisms are quasi-expanding
([9]). Quasi-expansion is the two dimensional analog of $f$ having
non-recurrent critical points. The exploitation of the properties
of quasi-expanding diffeomorphisms leads to the proofs of the
Theorems 5.1 and 5.2.

We believe that the connections made so far do not represent the
end of the story but only the beginning.  We trust that the
picture of two dimensional complex dynamics will become clearer
with time and, as it does, there will be valuable interactions
with the theory of real dynamics.

\head{References}

\item {[1]} E. Bedford \& J. Smillie, {\it Polynomial diffeomorphisms of
${\bf C}^2$: currents, equilibrium measures and hyperbolicity},
Inventiones Math. {\bf 103} (1991), 69--99.

\item {[2]} E. Bedford \& J. Smillie, {\it Polynomial diffeomorphisms of
${\bf C}^2$ II: stable manifolds and recurrence}, 4 No. 4 Journal
of the A.M.S. (1991), 657--679.

\item {[3]} E. Bedford \& J. Smillie, {\it Polynomial diffeomorphisms of
${\bf C}^2$ III: ergodicity, exponents and entropy of the
equilibrium measure}, Math. Annalen {\bf 294} (1992), 395--420.

\item{[4]} E.\ Bedford, M.\ Lyubich, and J.\ Smillie,  Polynomial
diffeomorphisms of ${\bf C}^2$ IV: the measure of maximal entropy
and laminar currents, {\it Inventiones Math.} 112 (1993), 77--125.

\item{[5]} E.\ Bedford, M.\ Lyubich, and J.\ Smillie,  Distribution of
periodic points of polynomial diffeomorphisms of ${\bf C}^2$,
Inventiones Math. {\bf 114} (1993), 277--288.

\item{[6]} E. Bedford \& J. Smillie, {\it Polynomial diffeomorphisms of ${\bf
C}^2$. V: Critical points and Lyapunov exponents},  J. Geom. Anal.
8 no. 3, (1998), 349--383.

\item{[7]} E. Bedford \& J. Smillie,  Polynomial diffeomorphisms
of ${\bf C}^2$. VI: Connectivity of $J$, {\it Annals of
Mathematics}, 148 (1998), 695--735.

\item{[8]} E. Bedford \& J. Smillie,  Polynomial diffeomorphisms
of ${\bf C}^2$. VII: Hyperbolicity and external rays, {\it Ann.\
Sci.\ Ecole Norm.\ Sup.\ 4 s\'erie}  32 (1999), 455--497.

\item{[9]}E. Bedford \& J. Smillie,  Polynomial diffeomorphisms
of ${\bf C}^2$. VIII: Quasi-expansion, {\it American Journal of
Mathematics}
  124 (2002), 221--271.

\item{[10]} E. Bedford \& J. Smillie,  Real polynomial diffeomorphisms with
maximal entropy: tangencies, (available at http://www.arXiv.org).

\item{[11]} M. Benedicks \& L. Carleson, The dynamics of the H\'enon map, {\it
Annals of Mathematics},  133, (1991), 73--179.

\item{[12]} H. Brolin, Invariant sets under iteration of rational functions,
{\it Ark. Mat}, 6, (1965), 103--144.

\item{[13]} A. de Carvalho \& T. Hall, How to prune a horseshoe, {\it
Nonlinearity}, 15 no. 3, (2002), R19--R68.

\item{[14]} R. Devaney \& Z. Nitecki,  Shift automorphisms in the H\'enon
mapping. Comm. Math. Phys. 67 (1979), no. 2, 137--146.

\item{[15]} S. Friedland \& J. Milnor, Dynamical properties of plane
polynomial automorphisms, {\it Ergodic Theory Dyn. Syst.} 9,
(1989), 67--99 .

\item{[16]} J. Hubbard \& R. Oberste-Vorth, H\'enon mappings in the complex
domain. II.  Projective and inductive limits of polynomials, in
{\it Real and complex dynamical systems} Kluwer, 1995.

\item{[17]} M. Lyubich, The quadratic family as a qualitatively solvable model
of chaos.  {\it Notices Amer. Math. Soc.} 47 (2000), no. 9,
1042--1052.

\item{[18]} A. Manning, The dimension of the maximal measure for a polynomial
map, {\it Ann. Math.} 119 (1984), 425--430.

\item{[19]} S. Newhouse, The abundance of wild hyperbolic sets and nonsmooth
stable sets for diffeomorphisms. {\it Inst. Hautes ¨¦tudes Sci.
Publ. Math.} No. 50 (1979), 101--151.

\item{[20]} F. Przytycki, Hausdorff dimension of the harmonic measure on the
boundary of an attractive basin for a holomorphic map, {\it
Invent. math.} 80 (1985), 161--179.

\end